\documentclass[titlepage,11pt]{article}
\usepackage{latexsym}
\usepackage{amsfonts}
\usepackage{amsmath}
\newcommand\blackslug{\hbox{\hskip 1pt \vrule width 4pt height 8pt depth 1.5pt
        \hskip 1pt}}
\newcommand\bbox{\hfill \quad \blackslug \bigbreak}
\def\d{\hbox{-}}
\def\cc{\hbox{-}\cdots\hbox{-}}
\def\l{,\ldots,}

\title{Induced subgraphs of graphs with large chromatic number.
\\VII. Gy\'arf\'as' complementation conjecture}
\author{Alex Scott\\
Mathematical Institute, University of Oxford, Oxford OX2 6GG, UK
\\
\\
Paul Seymour\thanks{Supported by ONR grant N00014-14-1-0084 and NSF
grant DMS-1265563.}\\
Princeton University, Princeton, NJ 08544, USA}

\date{September 23, 2016; revised \today}

\newtheorem{thm}{}[section]

\newcommand{\Proof}{\noindent{\bf Proof.}\ \ }

\begin{document}
\maketitle
\begin{abstract}
A class of graphs is {\em $\chi$-bounded} if there is a function $f$ such that $\chi(G)\le f(\omega(G))$ for every induced subgraph $G$ of every graph
in the class, where $\chi,\omega$ denote the chromatic number and clique number of $G$ respectively. 
In 1987, Gy\'arf\'as conjectured
that for every $c$, if $\mathcal{C}$ is a class of graphs such that $\chi(G)\le \omega(G)+c$ for every induced subgraph $G$
of every graph in the class, then the class of complements of members of $\mathcal{C}$ is $\chi$-bounded. We prove this conjecture. 
Indeed, more generally, a class of graphs is $\chi$-bounded if it has the property that no graph in the class has $c+1$
odd holes, pairwise disjoint and with no edges between them.
The main tool is a lemma that if $C$ is a shortest odd hole in a graph, and $X$ is the set of vertices 
with at least five neighbours in $V(C)$, then there is a three-vertex set that dominates $X$.
\end{abstract}

\section{Introduction}
All graphs considered are finite and simple. The chromatic number of $G$ is denoted by $\chi(G)$ and the cardinality of the largest
clique by $\omega(G)$. A class $\mathcal{C}$ of graphs is {\em $\chi$-bounded} if there is a function $f$ such that 
$\chi(H)\le f(\omega(H))$ for all $G\in \mathcal{C}$ and every induced subgraph $H$ of $G$. 
A graph $G$ is {\em perfect} if $\chi(H)=\omega(H)$ for every induced subgraph $H$ of $G$; and 
the class of perfect graphs is the primordial $\chi$-bounded class. Lov\'asz~\cite{lovasz} proved that the 
complement of every perfect graph is perfect, and 
it is natural to ask whether an analogue of Lov\'asz's theorem holds for other $\chi$-bounded classes.
Not always; for instance the class of all triangle-free graphs is not $\chi$-bounded, yet the class of their complements
is. 

Let $f$ be a function from the set of nonnegative integers into itself.
We say that $f$ is a {\em bounding function} for a graph $G$ if $\chi(H)\le f(\omega(H))$ for every
induced subgraph $H$ of $G$. Let $\mathcal{C}_f$ be the class of graphs $G$ such that $f$ is a bounding function for $G$.
We say $f$ has a {\em complementary bounding function} if the class of complements of members of $\mathcal{C}_f$
is $\chi$-bounded. Thus, Lov\'asz's theorem implies that the function $f(x)=x$ has a complementary bounding function, 
and we might ask which other functions do.

A. Gy\'arf\'as~\cite{gyarfas} proved that if $f$ has a complementary bounding function then $\inf f(x)/x=1$, and proposed
as a conjecture (conjecture 6.3 of~\cite{gyarfas}) 
that for all $c$ the function $f(x)=x+c$ has a complementary bounding function. That is our main result:

\begin{thm}\label{mainthm}
Let $c\ge 0$, and let $\mathcal{C}$ be a class of graphs such that 
$\chi(H)\le \omega(H)+c$ for all $G\in \mathcal{C}$ and every induced subgraph $H$ of $G$. Let $\mathcal{C}'$ be the class of 
complements of members of $\mathcal{C}$. Then $\mathcal{C}'$ is $\chi$-bounded.
\end{thm}
This question was already studied, for instance by Gy\'arf\'as, Li, Machado, Seb\H{o}, Thomass\'e and Trotignon,
in~\cite{gyarfas2}. They proved, among other things, that if $f$ is a function such that $f(x)=x$ for all sufficiently large $x$,
then $f$ has a complementary bounding function. (They also replaced the term ``binding function'' of \cite{gyarfas}
with ``bounding function'', which seems more appropriate, and we are following their terminology.)

Incidentally, one might wonder whether any classes with bounding function $f(x)=x+c$ exist; but there are several of interest. For instance:
\begin{itemize}
\item perfect graphs have bounding function $f(x)=x$; 
\item line graphs have a bounding function $f(x)=x+1$, by Vizing's theorem;
\item the intersection graph of a set of Jordan curves in the plane, pairwise intersecting in at most one point, has a bounding 
function $f(x)=x+327$, by a theorem of Cames van Batenburg, Esperet and M\"uller~\cite{esperet}; and
\item Esperet, Gon\c{c}alves and Labourel~\cite{esperet2} conjecture that for some $c$, $f(x)=x+c$ is a 
bounding function for 
intersection graphs of sets of non-crossing curves in the plane, pairwise intersecting in 
at most one point. 
\end{itemize}

A {\em hole} in a graph $G$ is an induced cycle of length at least four, and an {\em antihole} in $G$ 
is a hole in the complement graph $\overline{G}$. Two induced subgraphs of $G$ are {\em complete} 
if they are vertex-disjoint and every vertex of the first is adjacent in $G$ to every vertex of the second; and they are 
{\em anticomplete} if they are disjoint and there are 
no edges of $G$ joining them. If $G$ is a graph such that 
$\chi(H)\le \omega(H)+c$ for every induced subgraph $H$ of $G$, then there do not exist $c+1$ odd antiholes in $G$, 
pairwise complete.
Thus the following (applied to $\mathcal{C}'$) implies \ref{mainthm}.

\begin{thm}\label{mainthm2}
Let $c\ge 0$, and let $\mathcal{C}$ be a class of graphs $G$ such that there do not exist $c+1$ odd holes in $G$, pairwise anticomplete.
Then $\mathcal{C}$ is $\chi$-bounded.
\end{thm}

Vaidy Sivaraman points out (private communication) that \ref{mainthm2} also 
implies another conjecture of Gy\'arf\'as (conjecture 6.8 of~\cite{gyarfas}), 
the following, where $\alpha(H)$ denotes the cardinality of the largest 
stable subset of $H$:
\begin{thm}\label{gyarfasconj2}
Let $\mathcal{C}$ be the class of graphs $G$ such that 
$\alpha(H)\omega(H)\ge |V(H)|-1$ for every induced subgraph $H$ of $G$. 
Then $\mathcal{C}$ is $\chi$-bounded.
\end{thm}
This is immediate from \ref{mainthm2} since if $H$ is the disjoint union of 
two odd cycles of length at least four, then
$\alpha(H)\omega(H)=|V(H)| -2$.
  
Let $C$ be an odd hole in $G$. 
If $C$ has minimum length among all odd holes in $G$ we say $C$ is a {\em shortest} odd hole.
We denote the set of vertices in $V(G)\setminus V(C)$ with at least one neighbour in $V(C)$
by $N(C)$. If $X\subseteq V(G)$, we sometimes speak of $\chi(X)$ meaning $\chi(G[X])$.
As we shall see, the next result implies \ref{mainthm2}:

\begin{thm}\label{oddholenbr}
Let $G$ be a graph  such that for every vertex, its set of neighbours has chromatic number at most~$\tau$;
and let $C$ be a shortest odd hole in $G$. Then $N(C)$ has chromatic
number at most $21\tau$.
\end{thm}

We will also need the main result of~\cite{oddholes}, that:

\begin{thm}\label{oddhole}
Let $G$ be a graph with no odd hole.  Then $\chi(G)\le 2^{2^{\omega(G)+2}}$.
\end{thm}

\noindent{\bf Proof of \ref{mainthm2}, assuming \ref{oddholenbr}.\ \ }
Define $\tau(\kappa,c) = \kappa$ for $\kappa=0,1$. For $\kappa\ge 2$ let 
$\tau(\kappa,0) = 2^{2^{\kappa+2}}$, and let $\tau(\kappa,c)=\tau(\kappa,c-1)+21\tau(\kappa-1,c)$ for $c>0$.
We prove by induction on $\kappa+c$ that 
$\chi(G)\le \tau(\kappa,c)$ for every graph $G$ with $\omega(G)\le \kappa$ and
with no $c+1$ odd holes, pairwise anticomplete. This is true if $\kappa\le 1$, so we may assume that
$\kappa\ge 2$.
The claim is also true if $c=0$, by \ref{oddhole}, so we may assume that $c>0$.
Let $G$ be a graph with $\omega(G)\le \kappa$ and
with no $c+1$ odd holes, pairwise anticomplete, and suppose that $\chi(G)>\tau(\kappa,c)$. 
For every vertex $v$
of $G$, the set of neighbours of $v$ in $G$ has clique number at most $\kappa-1$, and hence 
has chromatic number at most $\tau(\kappa-1,c)$.
Since
$\tau(\kappa,c)\ge \tau(\kappa,0)= 2^{2^{\kappa+2}}$, there is an odd hole in $G$ by \ref{oddhole}; 
choose a shortest odd hole $C$. Let $G'$
be the graph obtained from $G$ by deleting $N(C)$. 
By \ref{oddholenbr}, $N(C)$
has chromatic number at most $21\tau(\kappa-1,c)$; and so 
$$\chi(G')\ge \chi(G)-21\tau(\kappa-1,c)>\tau(\kappa,c)-21\tau(\kappa-1,c) = \tau(\kappa,c-1).$$
Since $\chi(G')> 3 = \chi(C)$ and $C$ is a component of $G'$, it follows that $G'\setminus V(C)$ also has chromatic number more than
$\tau(\kappa,c-1)$.
But $G'\setminus V(C)$ does not contain $c$ pairwise anticomplete odd holes, a contradiction. 
This proves \ref{mainthm2}.~\bbox

Thus, it remains to prove \ref{oddholenbr}, and we will tackle this in the next two sections.
Another way prove \ref{mainthm2} is to prove directly the following strengthening of \ref{oddhole} (which also follows
from \ref{oddholenbr}):

\begin{thm}\label{oddholedel}
For all $\kappa, \tau\ge 0$ there exists $\tau'\ge 0$ such that if $G$ is a graph with $\omega(G)\le \kappa$ and
$\chi(G)>\tau'$, then there is an odd hole $C$ in $G$ such that if $X$ denotes the set of vertices of $G$ that are not in $C$ and have no neighbour in $C$, then $\chi(X)>\tau$.
\end{thm}
This can be proved by modifying slightly the proof of \ref{oddhole}; we check that
throughout the proof of \ref{oddhole}, whenever we find an odd hole, it can be chosen such that the set of vertices
with no neighbours in it still has large chromatic number.
It follows that in every graph $G$ with $\omega(G)\le \kappa$ and $\chi(G)$ sufficiently large, there exist
many odd holes pairwise anticomplete. To see this, choose an odd hole by \ref{oddholedel}, such that when we delete
it and its neighbours the chromatic number stays large, and repeat many times. And so \ref{mainthm2} would follow.

However, the proof via \ref{oddholenbr} is numerically better, and \ref{oddholenbr} may be of some independent interest.
For instance, the difficult part of the proof of correctness of the algorithm to test whether a graph is Berge, given
in~\cite{bergealg}, involved a closely related question which we discuss later.

\section{$C$-minor vertices}

Let $C$ be a shortest odd hole in $G$. 
Let us say a vertex $v\in N(C)$ is {\em $C$-minor}
if its neighbours in $C$ all belong to some path of $C$ of length two, and {\em $C$-major} otherwise. 
In this section we handle the $C$-minor vertices.

Let $v$ be a $C$-minor vertex, and choose a minimal path of $C$ of length at most two containing all neighbours of $v$ in $C$.
This path might have length zero, one or two, and we call this length the {\em type} of $v$.
We need the following lemma (its proof is trivial since all vertices of $K$ have degree at most $2|I|$):
\begin{thm}\label{cyclo}
Let $I\subseteq \{1\l n\}$, and let 
$K$ be a graph with vertex set $\{1\l n\}$, such that for all $i,j\in \{1\l n\}$ with $i<j$, if $i,j$ are adjacent in $K$ then 
either $j-i\in I$ or $i+n-j \in I$. Then $\chi(K)\le 2|I|+1$.
\end{thm}
We deduce

\begin{thm}\label{type0}
Let $G$ be a graph such that for every vertex, its set of neighbours has chromatic number at most $\tau$, and 
let $C$ be a shortest odd hole in $G$. The set of $C$-minor vertices of type zero induces a subgraph with chromatic
number at most $5\tau$.
\end{thm}
\Proof
Let $C$ have vertices $c_1\l c_{n}$ in order; and let $X_i$ be the set of all $C$-minor vertices of type zero that are adjacent to $c_i$,
for $1\le i\le n$. It follows that $\chi(X_i)\le \tau$ for $1\le i\le n$. Let $1\le i<j\le n$, and suppose that $u\in X_i$
is adjacent in $G$ to $v\in X_j$, and $c_i,c_j$ are not adjacent. 
One of the paths of $C$ between $c_i,c_j$ has even length, and its union with the path
$c_i\d u\d v\d c_j$ is an induced odd hole $C'$; and since $C$ is an odd hole of minimum length, $C'$ has length 
at least that of $C$.
Consequently the odd path of $C$ between $c_i,c_j$ has length three (if $c_i,c_j$ are not adjacent) or one; 
and so either $j-i=1$, or $j-i=3$, or
$i+n-j=1$, or $i+n-j=3$. Since this holds for all $i,j$, \ref{cyclo} (with $I=\{1,3\}$) implies that the union of the sets $X_i$
has chromatic number at most $5\tau$. This proves \ref{type0}.~\bbox

\begin{thm}\label{type1}
Let $G$ be a graph such that for every vertex, its set of neighbours has chromatic number at most $\tau$, and
let $C$ be a shortest odd hole in $G$. The set of $C$-minor vertices of type one induces a subgraph with chromatic
number at most $5\tau$.
\end{thm}
\Proof
Let $C$ have vertices $c_1\l c_{n}$ in order; and let $X_i$ be the set of all $C$-minor vertices of type one that are adjacent to $c_i, c_{i+1}$,
for $1\le i\le n$, where $c_{n+1}$ means $c_1$. It follows that $\chi(X_i)\le \tau$ for $1\le i\le n$. 
Let $1\le i<j\le n$, and suppose that $u\in X_i$
is adjacent in $G$ to $v\in X_j$. Let $P$ be the path with vertices $c_{i+1}\cc c_j$, and let $Q$ be the path with vertices
$c_{j+1}\cc c_n\d c_1\cc c_i$. One of $P,Q$ has odd length and one is even. If $Q$ is odd 
then 
the union of $P$ with the path
$c_{i+1}\d u\d v\d c_{j}$ is an induced odd cycle $C'$. Thus either $C'$ has length three or it has length at least that of $C$.
In the first case $j-i=1$, and in the second $Q$ has length one, that is, $i+n-j=2$. Similarly if $P$ is odd then either
$j-i=2$ or $i+n-j=1$. By \ref{cyclo} (with $I=\{1,2\}$)
the result follows.~\bbox

\begin{thm}\label{type2}
Let $G$ be a graph such that for every vertex, its set of neighbours has chromatic number at most $\tau$, and
let $C$ be a shortest odd hole in $G$. The set of $C$-minor vertices of type two induces a subgraph with chromatic
number at most $5\tau$.
\end{thm}
\Proof
Let $C$ have vertices $c_1\l c_{n}$ in order; and for $1\le i\le n$ 
let $X_i$ be the set of all $C$-minor vertices of type two that are adjacent to $c_{i-1}$ and to $c_{i+1}$ (and possibly to $c_i$), 
where $c_{n+1}$ means $c_1$ and $c_{-1}$ means $c_{n-1}$. It follows that $\chi(X_i)\le \tau$ for $1\le i\le n$.
Let $1\le i<j\le n$, and suppose that $u\in X_i$
is adjacent in $G$ to $v\in X_j$. It might be that $j=i+1$ or $i+n=j+1$; suppose not. 
Let $P$ be the path with vertices $c_{i+1}\cc c_{j-1}$, and let $Q$ be the path with vertices
$c_{j+1}\cc c_n\d c_1\cc c_{i-1}$. One of $P,Q$ has odd length and one is even. If $Q$ is odd
then
the union of $P$ with the path
$c_{i+1}\d u\d v\d c_{j-1}$ is an induced odd cycle $C'$, and $C'$ has length less than $C$. Thus $C'$ has length three, 
that is, $j-i=2$. Similarly if $P$ is odd then 
$i+n-j=2$. Thus we have shown that either $j-i=1$, or $i+n-j=1$, or $j-i=2$, or $i+n-j=2$.
By \ref{cyclo} (with $I=\{1,2\}$), the result follows.~\bbox

Adding, we deduce from \ref{type0}, \ref{type1} and \ref{type2} that:

\begin{thm}\label{minor}
Let $G$ be a graph such that for every vertex, its set of neighbours has chromatic number at most $\tau$, and
let $C$ be a shortest odd hole in $G$. The set of $C$-minor vertices induces a subgraph with chromatic
number at most $15\tau$.
\end{thm}

\section{$C$-major vertices}

If $u,v$ are $C$-major vertices, we say they form a {\em skew pair} if 
they are nonadjacent, each has exactly three
neighbours in $C$, and there is an odd path of $C$ with vertices $c_1\l c_n$ in order, where $n\ge 6$,
such that $u$ is adjacent to $c_1,c_2,c_{n-2}$,
and $v$ to $c_3,c_{n-1}, c_n$  (and hence $u,v$ have no common neighbour in $C$). 
If $X\subseteq V(G)$, we say that $Y$ {\em dominates} $X$ if every vertex in $X$
either belongs to $Y$ or has a neighbour in $Y$. We will prove:

\begin{thm}\label{hitmajor}
Let $C$ be a shortest odd hole in a graph $G$, and let $X$ be the set of $C$-major vertices.
\begin{itemize}
\item If some pair of vertices in $X$ is a skew pair, then there is a set of six vertices that dominates $X$.
\item If there is no skew pair in $X$, but some vertex in $X$ has at most four neighbours in $V(C)$, then there is 
a set of five vertices that dominates $X$.
\item If every vertex in $X$ has at least five neighbours in $V(C)$, then there is a set of three vertices that
dominates $X$.
\end{itemize}
\end{thm}

\noindent{\bf Proof of \ref{oddholenbr}, assuming \ref{hitmajor}.\ \ }
Let $G,\tau$ be as in \ref{oddholenbr}, and let $X$ be the set of $C$-major vertices.
By \ref{minor}, the set of $C$-minor vertices has chromatic number at most $15\tau$; and by \ref{hitmajor},
the set of $C$-major vertices has chromatic number at most $6\tau$. Consequently, $N(C)$ has chromatic number at
most $21\tau$, as required.~\bbox

Some of the results of this section are related to results of~\cite{bergealg}. In particular, one of the
goals of \cite{bergealg} was to find a small set of vertices which together dominated all $C$-major vertices, where
$C$ is a shortest odd hole. There are some differences: for instance,
in \cite{bergealg}, it was sometimes assumed that the graph contained none of a list 
of ``easily detectable subgraphs'', which we cannot, and in \cite{bergealg} the goal was to give (in polynomial time)
a list of subsets of $V(G)$ guaranteed to include the set of $C$-major vertices (without knowledge of $C$),
which we do not need. Nevertheless,
some of the proofs here are derived from proofs in~\cite{bergealg}.

Let $C$ be a cycle, and let $A\subseteq V(C)$. An {\em $A$-gap} is a subgraph of $C$ composed of a component $X$
of $C\setminus A$, the vertices of $A$ with neighbours in $X$, and the edges between $A$ and $X$.
(So if some two vertices in $A$ are nonadjacent, the $A$-gaps are the paths of $C$ of length $\ge 2$, with both ends
in $A$ and no internal vertex in $A$.) If $A,B\subseteq V(C)$, an {\em $(A,B)$-gap} is a path $P$ of $C$ between $a,b$ say, such that $a$
is the unique vertex of $P$ in $A$ and $b$ is the unique vertex of $P$ in $B$. (Possibly $a = b$ and
$P$ has length $0$.)
The {\em length} of an $A$-gap or $(A,B)$-gap is the number of edges in it (so if
$A$ consists just of two adjacent vertices, there is a single $A$-gap, of length $|E(C)| - 1$).
Every $A$-gap has length at least two, but an $(A,B)$-gap can have length $0$ or $1$ as well.
We speak of cycles, $A$-gaps, and $(A,B)$-gaps being {\em odd} or {\em even} meaning that
they have an odd (or even, respectively) number of edges.
We say that $A$ is {\em normal} in $C$ if every $A$-gap is even. (When $A=\emptyset$, $C$ is an $A$-gap, and therefore
$\emptyset$ is not normal.) An {\em $A$-edge} is an edge of $C$ with both ends in $A$.
Thus every edge of $C$ is either an $A$-edge or belongs to a unique $A$-gap (and not both).

To prove \ref{hitmajor} we first need the following.

\begin{thm}\label{getgap}
Let $C$ be an odd cycle, and let $A,B\subseteq V(C)$ be normal. Then 
every odd $(A,B)$-gap is disjoint from every even $(A,B)$-gap. Moreover, if $A\cap B$ is not normal, then
there is an odd $(A,B)$-gap and an even $(A,B)$-gap.
\end{thm}
\Proof
Let $P,Q$ be odd and even $(A,B)$-gaps respectively, and suppose that there exists $v\in V(P\cap Q)$.
Then $v$ is an end of both $P,Q$, and belongs to one of $A,B$, say $A$; but
then $P\cup Q$ is an
odd $B$-gap, a contradiction. This proves the first assertion of the theorem.

Now we assume that $A\cap B$ is not normal, and so there is an odd $A\cap B$-gap $R$ say.
Every edge of $R$ is either an $A$-edge or belongs to an $A$-gap, and since $R$ is odd and all $A$-gaps are even it
follows that an  odd number of edges of $R$ are $A$-edges. Since every $A$-edge of $R$ belongs to a $B$-gap,
some $B$-gap $P$ included in $R$ contains an odd number of $A$-edges. Since $P$ is even, there is an odd number of edges
of $P$ that are not $A$-edges. Each such edge belongs either to a $A$-gap included in $P$, or to an $(A,B)$-gap included
in $P$.
The number of edges of $P$ that belong to $A$-gaps included in $P$ is even, since all $A$-gaps are even; so an odd number
of edges of $P$ belong to $(A,B)$-gaps included in $P$. Since there are exactly two $(A,B)$-gaps included in $P$, one of
them is odd and one is even. 
This proves \ref{getgap}.~\bbox

\begin{thm}\label{abnormal2}
Let $C$ be a shortest odd hole in $G$, and let $u,v$ be nonadjacent $C$-major vertices.
Let $A,B$ be the sets of neighbours in $V(C)$ of $u,v$ respectively. Then  either:
\begin{itemize}
\item $A\cap B$ is normal; or
\item $u,v$ form a skew pair; or
\item $A\cap B\ne\emptyset$ and $\min(|A|,|B|)=3$; or
\item there is a path $c_1\cc c_5$ of $C$ such that $\{c_1\l c_4\}$ is the set of neighbours of $u$ in $V(C)$, and $\{c_2\l c_5\}$
is the set of neighbours of $v$ in $V(C)$.
\end{itemize}
\end{thm}
\Proof
Since $u,v$ are $C$-major and $C$ is a shortest odd hole it follows that $A,B$ are normal.
We will show that we may assume that:
\\
\\
(1) {\em If $P$ is an odd $(A,B)$-gap and $Q$ is an even $(A,B)$-gap then $V(P\cap Q)=\emptyset$, and there is an
edge $e$ of $C$ joining an end of $P$ with an end of $Q$, and $e$ is either an $A$-edge or a $B$-edge.}
\\
\\
By \ref{getgap}, $V(P\cap Q)=\emptyset$.
Suppose first that $P,Q$ are anticomplete. Then the subgraph induced on $V(P\cup Q)\cup\{u,v\}$ is an odd hole, and so its length
is at least that of $C$; and hence at most two vertices of $C$ do not belong to $P\cup Q$. 
Thus we can number the vertices of $C$ in order as
$c_1\l c_n$ where $P$ is the subpath $c_1\d c_2\cc c_r$, and $Q$ is the subpath $c_{r+2}\cc \ c_{n-1}$ (where $2\le r\le n-3$,
and possibly $r=n-3$).
Since $u,v$ both have only
two neighbours in $P\cup Q$ it follows that $|A|,|B|\le 4$. First suppose that not both $c_{r+1},c_n$ belong to $A\cap B$.
Hence $\min(|A|,|B|)=3$, and so we may assume that $A\cap B=\emptyset$, and therefore
$u,v$ are each adjacent to exactly one of $c_{r+1}, c_n$.
From the symmetry we may assume that $c_{n-1}\in A$, and hence $c_{r+2}\in B$. Consequently $c_{r+1}\notin A$ 
since $c_{r+1}\cc c_{n-1}$ is not an odd $A$-gap. Thus $c_{r+1}\in B$,
and $c_n\in A$.  Since
$u$ is $C$-major, $A\not\subseteq \{c_{n-1}, c_n, c_1\}$; and so $c_r\in A$ and $c_1\in B$, and
so $u,v$ form a skew pair (the path $c_{n-1}\d c_n\d c_1\cc c_r\d c_{r+1}\d c_{r+2}$ corresponds to the path
numbered $c_1\cc c_n$ in the definition of a skew pair).

We may therefore assume that $c_{r+1}, c_n\in A\cap B$.
By exchanging $A,B$ if necessary, we may assume that $c_{r+2}\in B$, and since $c_{r+2}\cc c_n$ is not an odd $B$-gap, it follows that
$r+2=n-1$, and so $c_{n-1}\in A\cap B$. But then the last outcome of the theorem holds.

Thus we may assume that $P,Q$ are not anticomplete.
Hence some end $p$ of $P$ is adjacent to some end $q$ of $Q$.
If one of
$p,q$ is in $A\setminus B$ and the other in $B\setminus A$, then $Q$ has positive length and
the subgraph induced on $V(Q)\cup \{p\}$ is either an odd $A$-gap or an odd $B$-gap,
in either case a contradiction. So one of $A,B$ contains both of $p,q$. This proves (1).

\bigskip

We may assume that $A\cap B$ is not normal. By \ref{getgap},
there is an odd $(A,B)$-gap $P$ and an even $(A,B)$-gap $Q$.
From (1), we may number
the vertices of $C$ as $c_1\l c_n$ in order such that $P$ is $c_1\cc  c_i$, and $Q$ is $c_{i+1}\cc c_j$, 
and one of $A,B$ contains both of $c_i,c_{i+1}$. We may assume, exchanging $A,B$ if necessary, 
that $c_1, c_j\in A$ and $c_i, c_{i+1}\in B$. (Possibly $i+1=j$.)
Now one of $c_{j+2}\l c_n\in B$; for otherwise $c_{j+1}\in B$ (since $|B|\ge 3$), and $i+2\ne j$ (since $v$ is $C$-major),
and $c_{i+1}\cc c_{j+1}$ is an odd $B$-gap, which is impossible. 
Since $c_1\in A\setminus B$, there is an $(A,B)$-gap $S$ (possibly of length zero) that is a subpath of the path $c_{j+2}\cc c_n\d c_1$.
Now $S$ is anticomplete to $Q$, since $c_i,c_{j+1}\notin V(S)$; so by (1) applied to $S,Q$, 
it follows that $S$ is even. By (1) applied to $P,S$, $c_1\notin V(S)$, and $c_n$ is an end of $S$, and $c_n\in A$.

Suppose that both $A,B$ have nonempty intersection with $\{c_{i+2}\l c_{n-1}\}$. Then there is an $(A,B)$-gap $R$ say
with vertex set a subset of this set. Since $R$ is anticomplete to $P$, (1) implies that $R$ is odd. By (1), $R$ is disjoint
from both $Q,S$, and the union of the vertex sets of $P,Q,R,S$ equals $V(C)$, which is impossible since $C$ is odd.
Thus not both $A,B$ have nonempty intersection with $\{c_{i+2}\l c_{n-1}\}$. This implies that $\min(|A|,|B|)=3$, and also that
one of $Q,S$
has length zero, so $A\cap B\ne \emptyset$ and the theorem holds. This proves \ref{abnormal2}.~\bbox

We deduce:
\begin{thm}\label{commonnbr}
Let $C$ be a shortest odd hole in a graph $G$, and let $u,v$ be $C$-major vertices. Then either $u,v$
are adjacent, or they form a skew pair, or some vertex in $C$ is adjacent to both $u,v$.
\end{thm}
\Proof
Let $A,B$ be the sets of neighbours in $V(C)$ of $u,v$ respectively. Thus, $A,B$ are both normal, since $C$ is a shortest odd hole
and $u,v$ are $C$-major. Suppose that $u,v$ are nonadjacent and $A\cap B=\emptyset$ and $u,v$ do not form a skew pair.
In particular, $A\cap B$ is not normal, contrary to \ref{abnormal2}. This proves \ref{commonnbr}.~\bbox

Next we need the following
(compare theorem 7.6 of \cite{bergealg}):

\begin{thm}\label{stablenormal}
Let $C$ be a shortest odd hole in $G$, and let $u,v$ be nonadjacent $C$-major vertices, both 
with at least five neighbours 
in $V(C)$. Let $A,B$ be the sets of neighbours of $u,v$ in $V(C)$ respectively. Then every $(A,B)$-gap is even.
\end{thm}
\Proof
Suppose that $P$ is an odd $(A,B)$-gap. Let
$C$ have vertices $c_1\l c_{n}$ in order, where $P$ is $c_1\cc c_r$, and $r$ is even, with $2\le r <n$.
\\
\\
(1) {\em If $Q$ is an even $(A,B)$-gap, then $V(P\cap Q) = \emptyset$ and there is an edge between $V(P),V(Q)$.}
\\
\\
For $P\cap Q$ is empty, by \ref{getgap},
so we assume there are no edges between $P$ and $Q$. Then their union, together with $u,v$, forms an
odd hole $C'$. Since $|A|\ge 5$, there are at least three vertices of $C$ that do not belong
to $P\cup Q$, and so $C'$ is shorter than $C$, a contradiction.
This proves (1).
\\
\\
(2) {\em There is an even $(A,B)$-gap.}
\\
\\
For if $A\cap B\ne \emptyset$ then there is
an $(A,B)$-gap of length $0$, so we assume that $A\cap B=\emptyset$. In particular $A\cap B$ is not normal, 
and the claim follows from \ref{getgap}. 
This proves (2).

\bigskip

By (1) and (2) we may assume that $c_{r+1}\cc c_s$ is an even $(A,B)$-gap, $Q$ say, for some odd $s$ with 
$r+1 \le s \le n$.
Since $|A|, |B|\ge 5$, and only two vertices of $P \cup Q$ belong to $A$, and the same for $B$, it follows that $s \le n-3$,
and both $A,B$ meet the path $c_{s+2}\cc c_{n-1}$, and so 
there is an $(A,B)$-gap contained in this path.
It is not even, by (1), since there are no edges between it and $P$, and it is not odd, by (1), since there
are no edges between it and $Q$, a contradiction. This proves \ref{stablenormal}.\bbox

We need theorem 7.4 of \cite{bergealg}, the following:
\begin{thm}\label{stablegap}
Let $C$ be an odd cycle. Let $A_1\l A_k\subseteq V(G)$ be normal, such that for $1 \le i < j \le k$, every
$(A_i,A_j)$-gap is even. Then $A_1\cap \cdots \cap A_k$ is normal.
\end{thm}

By combining \ref{stablenormal} and \ref{stablegap} we deduce:
\begin{thm}\label{stablenbrs}
Let $C$ be a shortest odd hole in a graph $G$, and let $X$ be a stable set of $C$-major vertices, 
each with at least five neighbours in $V(C)$.
Then the set of vertices in $V(C)$ that are complete to $X$ is normal.
\end{thm}

\begin{thm}\label{antihole}
Let $C$ be a shortest odd hole in $G$, let $z$ be a $C$-major vertex, and let $X$ be a set of $C$-major vertices
all different from and nonadjacent to $z$, and such that each vertex in $X\cup \{z\}$ has at least five 
neighbours in $V(C)$. Suppose that the set of common neighbours of $X\cup \{z\}$ in $V(C)$ is
not normal. Then there exist adjacent $u,v\in X$, and a path $c_1\d c_2\d c_3\d c_4$ of $C$, such that every two vertices
in $\{c_1,c_2,c_3,c_4,u,v,z\}$ are adjacent except for the pairs $c_3c_1,c_1c_4,c_4c_2,c_2u,uz,zv,vc_3$, which are 
nonadjacent.
\end{thm}
\Proof
Since the set of common neighbours of $X\cup \{z\}$ in $V(C)$ is
not normal, at least one vertex in $X\cup \{z\}$ has a nonneighbour in $V(C)$, and since they each have at least
five neighbours in $C$ it follows that $C$ has length at least six.
We proceed by induction on $|X|$. Since the set of common neighbours of $X\cup \{z\}$ in $V(C)$ is not normal,
\ref{stablenbrs} implies that some two vertices in $X$ are adjacent, say $u,v$. Let $A,B$ be the sets of 
common neighbours
in $V(C)$ of $(X\setminus \{u\})\cup \{z\}, (X\setminus \{v\})\cup \{z\}$ respectively. From the inductive hypothesis
we may assume that $A,B$ are normal; and $A\cap B$ is not normal. Consequently $A\setminus B$ and $B\setminus A$ are 
both nonempty. If $a\in A\setminus B$ and $b\in B\setminus A$
then $a$ is adjacent to $v$ and not to $u$, and $b$ is adjacent to $u$ and not to $v$; and since $z\d a\d v\d u\d b\d z$
is not a 5-hole (because $C$ is a shortest odd hole and has length at least six), it follows that $a,b$ are
adjacent. Hence every vertex in $A\setminus B$ is adjacent to every vertex in $B\setminus A$; and since
$A\setminus B$ and $B\setminus A$ are both nonempty, one of $A\setminus B,B\setminus A$ has cardinality one,
say $A\setminus B$, and so we may number the vertices of $C$ in order as $c_1\l c_n$, such that
$A\setminus B=\{c_2\}$, and $c_3\in B\setminus A\subseteq \{c_1,c_3\}$. Since $A$ is normal and hence $|A|\ge 2$, 
it follows that $A\cap B\ne \emptyset$. Since $A\cap B$ is not normal, there is an odd $A\cap B$-gap $R$ say. Since
$R$ is not an $A$-gap because $A$ is normal, its interior contains a member of $A\setminus B$, and since 
$A\setminus B=\{c_2\}$, we deduce that $c_2$ belongs to the interior of $R$. One of the paths of $R$ between $c_2$ and its ends is odd, and is therefore not an $A$-gap, and so has length one. Consequently one of $c_1,c_3$ is an end of $R$
and therefore belongs to $A\cap B$. Since $c_3\in B\setminus A$, it follows that $c_1\in A\cap B$, and so $B\setminus A=\{c_3\}$
By the same argument with $A,B$ exchanged, it follows that $c_4\in A\cap B$, and then the theorem is satisfied. This 
proves \ref{antihole}.~\bbox

\begin{thm}\label{hitmajor2}
Let $C$ be a shortest odd hole in $G$, let $z$ be a $C$-major vertex, and let $X$ be a set of $C$-major vertices
all different from and nonadjacent to $z$, and such that each vertex in $X\cup \{z\}$ has at least five
neighbours in $V(C)$. Then there is an edge $ab$ of $C$, such that $a,b$ are both adjacent to $z$, and every vertex
in $X$ is adjacent to at least one of $a,b$.
\end{thm}
\Proof
We proceed by induction on $|X|$. If the set of common neighbours of $X\cup \{z\}$ is normal then there is an
edge of $C$ with both ends complete to $X\cup \{z\}$, which therefore satisfies the theorem. Consequently we may assume 
that the set of common neighbours of $X\cup \{z\}$ is not normal. By \ref{antihole} there exist adjacent $u,v\in X$, and a path $c_1\d c_2\d c_3\d c_4$ of $C$, such that every two vertices
in $\{c_1,c_2,c_3,c_4,u,v,z\}$ are adjacent except for the pairs $c_3c_1,c_1c_4,c_4c_2,c_2u,uz,zv,vc_3$, which are
nonadjacent. From the inductive hypothesis there is an edge $ab$ of $C$, such that $a,b$ are both adjacent to $z$, and every vertex
in $X\setminus \{u\}$ is adjacent to at least one of $a,b$. If $u$ is also adjacent to one of $a,b$ then the theorem is
satisfied, so we assume not. Let $v$ be adjacent to $a$ say. Since $u$ is not adjacent to $a,b$ it follows that 
$a,b\ne c_1,c_3,c_4$; and therefore $a,b\ne c_2$ since $a,b$ are consecutive in $C$. But then $z\d c_3\d u\d v\d a\d z$
is a 5-hole, a contradiction. This proves \ref{hitmajor2}.~\bbox

We deduce \ref{hitmajor}, which we restate.

\begin{thm}\label{hitmajoragain}
Let $C$ be a shortest odd hole in a graph $G$, and let $X$ be the set of $C$-major vertices.
\begin{itemize}
\item If some pair of vertices in $X$ is a skew pair, then there is a set of six vertices that dominates $X$.
\item If there is no skew pair in $X$, but some vertex in $X$ has at most four neighbours in $V(C)$, then there is
a set of five vertices that dominates $X$.
\item If every vertex in $X$ has at least five neighbours in $V(C)$, then there is a set of three vertices that
dominates $X$.
\end{itemize}
\end{thm}
\Proof
Suppose first that some pair $u,v$ of vertices in $X$ is a skew pair. Consequently
each has exactly three
neighbours in $C$, and there is an odd path of $C$ with vertices $c_1\l c_n$ in order, where $n\ge 6$,
such that $u$ is adjacent to $c_1,c_2,c_{n-2}$,
and $v$ to $c_3,c_{n-1}, c_n$. We claim that $\{c_0,c_1,c_2,c_3,c_n,u\}$ dominates $X$ (where $c_0$ is the second vertex 
adjacent to $c_1$). 
For let $w\in  X\setminus \{u\}$.
If $w$ is adjacent to $u$ or has a common neighbour with $u$ then one of $c_1,c_2,c_{n-2},u$ is adjacent to $w$, and 
otherwise by \ref{commonnbr} $u,w$ form a skew pair, and so $w$ is adjacent to one of $c_{3},c_n$. Thus the result
holds if there is a skew pair. 

Now we assume that there is no skew pair, but some vertex $v$ in $X$ has at most four neighbours in $V(C)$. 
By \ref{commonnbr} the set consisting of $v$ and its neighbours in $V(C)$ dominates $X$, and again the theorem holds.
Thus we may assume that every vertex in $X$ has at least five neighbours in $V(C)$.

Choose $z\in X$, and let $X'$ be the set of vertices in $X\setminus \{z\}$ that are nonadjacent to $z$. By \ref{hitmajor2},
there is an edge $ab$ of $C$ such that every vertex in $X'$ is adjacent to one of $a,b$. But then $\{a,b,z\}$ dominates $X$. This proves \ref{hitmajoragain}.~\bbox

It is an entertaining problem to attempt to replace the 21 in \ref{oddholenbr} by as small a number as possible. There are 
a number of tricks that give savings, and here we sketch some of them, without full details.
\begin{itemize}
\item In \ref{hitmajoragain}, the six of the first bullet can be reduced to five; the set $\{c_1,c_2, c_3,c_{n-2},u\}$
dominates $X$ (because $u$ cannot be part of a skew pair on both sides, unless the shortest odd hole has length nine,
and in this case we can replace 21 by 9).
\item \ref{hitmajor2} can be extended to include the case when some $C$-major vertices have only four neighbours in $C$, provided
they are not all four consecutive.
\item In all cases of \ref{hitmajoragain}, there is a set of five vertices that dominates $X$ including three consecutive
vertices of $C$. To see this claim, we need to improve the second bullet of \ref{hitmajoragain}. 
When there is no skew pair, and some $C$-major vertex has
three neighbours in $C$, the claim is true since two of its neighbours are consecutive; otherwise, if some $v$ has only
four neighbours and they are consecutive in $C$, then $v$ plus the middle two of the four dominate $X$; and if there is no such $v$,
the conclusion of the third bullet applies (because of the strengthening of \ref{hitmajor2} mentioned above).
This will help with colouring the $C$-minor vertices.
\item Thus there are three consecutive vertices $c_{n-2}, c_{n-1}, c_n$ of $C$ as above, and we are going to pay $3\tau$ 
to colour the $C$-major 
vertices adjacent to one of the three. We might as well colour the $C$-minor vertices adjacent to one of $c_{n-2}, c_{n-1}, c_n$ at the same time;
and so we need not bother to colour those $C$-minor vertices in the course of \ref{type0}, \ref{type1}, \ref{type2}.
\item For \ref{type0}, we only need colour the union of the sets
$X_1\l X_{n-4}$. For $1\le i<j\le n-4$, if there is an edge between $X_i,X_j$ then $j-i\in \{1,3\}$ and in particular, $j-i$
is odd; so the union of $X_1\l X_{n-4}$ has chromatic number at most $2\tau$.
For \ref{type1}, we only need colour the union of the sets
$X_1\l X_{n-5}$. Each $X_i$ has neighbours in only two of the sets $X_{i+1}\l X_{n-5}$, and so this union has chromatic number
at most $3\tau$. Thus the $5\tau$ of \ref{type1} can be reduced to $3\tau$; and similarly the $5\tau$ of \ref{type2}
can be reduced to $3\tau$.
\item There is another improvement possible; it is more efficient to handle the $C$-minor vertices of type zero and type two 
simultaneously. Let $Y$ be the set of $C$-minor vertices of type zero or two that are nonadjacent to $c_{n-2},c_{n-1}, c_n$.
For $1\le i\le n-3$ let $X_i$ be the type zero vertices in $Y$ adjacent to $c_i$, together with the type two vertices in $Y$
adjacent to $c_i, c_{i+2}$. For $i<j$, if there is an edge between $X_i, X_j$ then $j\le i+3$; so the union of $X_1\l X_{n-3}$
has chromatic number at most $4\tau$. 
This shows we can colour all $C$-minor vertices of types zero or two nonadjacent to $c_{n-2}, c_{n-1}, c_n$ with $4\tau$ colours.
\end{itemize}
Thus we can colour all the $C$-minor vertices nonadjacent to $c_{n-2}, c_{n-1}, c_n$ with $7\tau$ colours. 
Handling the $C$-major vertices
and the vertices adjacent to $c_{n-2}, c_{n-1}$ or $c_n$ takes $5\tau$ more; so altogether we use $12\tau$ colours. Thus the 
21 of \ref{oddholenbr} can be reduced to 12.

\section{Acknowledgement}

We would like to thank Vaidy Sivaraman for pointing out Gy\'arf\'as' question and for some very helpful discussion; and Sang-Il Oum and his students for carefully reading the paper and finding several errors.

\end{document}